\documentclass[tikz]{article}

\usepackage{PRIMEarxiv}

\usepackage[utf8]{inputenc} 
\usepackage[T1]{fontenc}    
\usepackage{hyperref}       
\usepackage{url}            
\usepackage{booktabs}       
\usepackage{amsfonts}       
\usepackage{nicefrac}       
\usepackage{microtype}      
\usepackage{lipsum}
\usepackage{fancyhdr}       
\usepackage{graphicx}       

\usepackage{mathtools} 
\usepackage{booktabs} 
\usepackage{algorithm}
\usepackage[noend]{algorithmic}
\usepackage{multirow}
\usepackage{caption}
\usepackage{subcaption}
\usepackage{amssymb}
\usepackage{graphicx}
\usepackage{xfrac}
\usepackage{xcolor}
\usepackage{ulem}
\usepackage{microtype}
\title{
    Maximum Dispersion, Maximum Concentration: Enhancing the Quality of MOP Solutions
}

\author{
Gladston Moreira$^1$, Ivan Meneghini$^2$ and Elizabeth Wanner$^3$ \\
 $^1$Computing Department, Universidade Federal de Ouro Preto, Ouro
Preto, 35402-136, Minas Gerais, Brazil.\\
$^2$Federal Institute of Minas Gerais, Ibirité-MG, Brazil.\\
$^3$School of Computer Science and Digital Technologies, Aston University, UK.\\
    \texttt{gladston@ufop.edu.br, ivan.reinaldo@ifmg.edu.br, e.wanner1@aston.ac.uk}
}

\begin{document}
\maketitle

\begin{abstract}
    Multi-objective optimization problems (MOPs) often require a trade-off between conflicting objectives, maximizing diversity and convergence in the objective space. This study presents an approach to improve the quality of MOP solutions by optimizing the dispersion in the decision space and the convergence in a specific region of the objective space. Our approach defines a Region of Interest (ROI) based on a cone representing the decision maker's preferences in the objective space, while enhancing the dispersion of solutions in the decision space using a uniformity measure. Combining solution concentration in the objective space with dispersion in the decision space intensifies the search for Pareto-optimal solutions while increasing solution diversity. When combined, these characteristics improve the quality of solutions and avoid the bias caused by clustering solutions in a specific region of the decision space. Preliminary experiments suggest that this method enhances multi-objective optimization by generating solutions that effectively balance dispersion and concentration, thereby mitigating bias in the decision space.
\end{abstract}

\keywords{Multi-objective optimization \and Evolutionary algorithms \and Decision space diversity \and Region of Interest.}

\section{Introduction}\label{sec:intro}

Multi-objective optimization problems (MOPs) arise in various real-world applications where multiple conflicting objectives must be optimized simultaneously \cite{deb2001multiobjective}. These problems are common in fields such as engineering design \cite{Tawhid2018}, finance \cite{He2022}, logistics \cite{Luo2022}, neural networks \cite{Kaoutar2017} and natural science \cite{Wang2023}. Due to their ability to approximate the Pareto front efficiently, Multi-Objective Evolutionary Algorithms (MOEAs) have become widely used for solving MOPs \cite{mops2024}.

In any application formulated as a MOP, the decision maker seeks ``good solutions''. The presence of uncertainties in the decision space \cite{Meneghini2016} and/or noise in the objective space \cite{Sousa2024,Sousa2023} degrades the solutions obtained by traditional MOEAs. Such situations require specialized optimizers that yield effective solutions in uncertain scenarios. The studies in \cite{UD:24,KIDD2020,Jesus2018,Vaz:15,Eus:14} focus on approximating the Pareto set in bi-objective combinatorial optimization problems, evaluating the quality of the representation in the objective space based on uniformity, coverage, and the $\epsilon$-indicator measure. On the other hand, polarized solutions within a particular region of the decision space can cause bias in implementing a specific solution. The dispersion of efficient solutions in the decision space is evaluated in \cite{moreira:19,Narukawa:13,Ulr:10}. In \cite{Rob2022}, the authors define a robustness measure in decision space for multi-objective integer linear programming problems. 

A well-distributed set of solutions in the decision space enhances decision-making by offering a rich and diverse range of options. At the same time, prioritizing the search for solutions in a specific region of the objective space intensifies the selective pressure in Evolutionary Algorithms (EA) and improves convergence rates \cite{ivan1}. The combination of dispersion and concentration of solutions in the decision and objective spaces, respectively, is particularly productive when dealing with large-scale MOPs, where computational resources and the efficiency of MOEAs can become limiting factors.

Despite its importance, the combination of solution dispersion and concentration in MOEAs has not been extensively studied in the literature. Traditional performance metrics, such as hypervolume \cite{Shang2021} and inverse generational distance (IGD) \cite{CoelloCoello2004}, primarily focus on the convergence and coverage of the Pareto Front in objective space, but do not address solution diversity and representativeness in decision space. In practical applications, decision-makers typically seek a subset of optimal solutions within a specific region of the objective space, expecting to have diverse options for these solutions in the decision space.

This work introduces a new perspective on improving the quality of MOP solutions by maximizing dispersion and concentration, thereby balancing diversity in the decision space with accuracy in the objective space. We propose an adaptive mechanism that enhances the distribution and convergence of solutions in regions of interest within the objective space, while maintaining solution diversity in the decision space. By leveraging this approach, we aim to improve the robustness of MOEA-generated solutions in real-world decision-making scenarios.
%

\section{Background}

A multi-objective optimization problem (MOP) involves simultaneously optimizing multiple objective functions subject to constraints. In general, MOPs have conflicting objectives, meaning that improving one objective may lead to the deterioration of another. A MOP aims to find solutions that balance all objectives, as obtaining the optimal value for each objective is inherently impractical. In this study, optimization is approached as the minimization of objective values.

It is common to represent a MOP using the application \cite{deb2001multiobjective}
\begin{equation} \label{eq:mop1}
	\mathbf{F}: \Omega \subset \mathbb{R}^n \to \mathbb{R}^m
\end{equation}
In this representation,
$$\mathbf{x}=(x_1,x_2, \ldots, x_n) \in \Omega \subset \mathbb{R}^n$$
is a decision variable vector with $n$ components, $\mathbb{R}^n$ is called decision space and $\Omega$ is the feasible set of the decision space. Furthermore,
$$ \mathbf{F}(\mathbf{x}) = \left(f_1(\mathbf{x}), f_2(\mathbf{x}), \ldots, f_m(\mathbf{x})\right) $$
is a vector in $\mathbb{R}^m$, where $\mathbb{R}^m$ is called objective space and $\mathbf{F}(\Omega)$ is referred to as the feasible set image. Each of the $m$ coordinates $f_i(\mathbf{x})$ of the vector $\mathbf{F}(\mathbf{x})$ represents an objective to be optimized. Additionally, the application in \eqref{eq:mop1} can be subject to constraints in $\Omega$ and/or in $\mathbf{F}(\Omega)$ \cite{Liang2024}.

Since the space $\mathbb{R}^k$ with $k > 1$ lacks a total order relation that allows classification $\mathbf{a} \leq \mathbf{b}$ as is common for $k = 1$, a partial order relation $\prec$ is used to compare any two vectors in $\mathbb{R}^k$. The relation $\prec$ is called the Pareto dominance relation and is defined as follows: Let $\mathbf{a} = (a_1,a_2, \ldots, a_k)$ and $\mathbf{b} = (b_1, b_2, \ldots, b_k)$ be vectors in $\mathbb{R}^k$. If $a_i \leq b_i$ for $i = 1,2, \ldots , k$ and $\mathbf{a} \ne \mathbf{b}$, then $\mathbf{a} \prec \mathbf{b}$. In this case, $\mathbf{a}$ dominates $\mathbf{b}$. If $\nexists~\mathbf{c} \prec \mathbf{a}$, then $\mathbf{a}$ is non-dominated vector. In the context of a MOP, the order relation $\prec$ defines its set of solutions. To this end, consider the sets.
%
\begin{align}
	PF & = \{ \mathbf{y} \in F(\Omega); \nexists ~ \mathbf{x} \in \Omega \text{ such that  } F(\mathbf{x}) \prec \mathbf{y} \} \\
	PO & = \{ \mathbf{x} \in \Omega; F(\mathbf{x}) \in PF\}
\end{align}

The set $PF \subset F(\Omega)$ is called the Pareto Front, and $PO \subset \Omega$ is called the Pareto Optimal set or efficient solution set. These two sets define the solutions of a MOP, referred to as optimal or efficient solutions. Generally, $PF$ is a surface of dimension $m-1$ in $\mathbb{R}^m$.

Obtaining the sets $PF$ and $PO$ is a challenging task. The complexity of the application $\mathbf{F}$ and the magnitude of $n$ and $m$ in the application \eqref{eq:mop1} make deterministic methods (such as gradient descent and others) impractical. One way to obtain solutions for a MOP is to apply metaheuristics, particularly Multi-Objective Evolutionary Algorithms (MOEAs). These techniques do not guarantee the exact $PF$ and $PO$ sets but produce a finite set of non-dominated points (in the objective space) that are sufficiently close to the $PF$ set\footnote{The proximity between the set of solutions obtained by MOEA and the PF set is called convergence of solutions.}.

In addition to the computational complexity of a MOP, obtaining solutions via evolutionary algorithms also faces several challenges. In the decision space $\Omega$, the solutions obtained by an MOEA may cluster in specific regions, poorly representing the $PO$ set. Decision-makers who utilize the solutions obtained for the problem that generated the MOP expect to get a diverse, unbiased, and representative set of solutions in the decision space. These characteristics make the obtained solution robust, as others can replace less attractive solutions with different configurations in the decision space. A clustered set of solutions in the decision space lowers solution quality by limiting the available implementation options. If the initial conditions change after optimization, a non-diverse solution set may become infeasible, necessitating another round of MOEA execution. Due to its computational cost, this new search results in project delays and even financial costs. 

\subsection{Region of Interest}

Despite its importance, the issue of diversity and representativeness of solutions in the decision space is sparsely addressed in the literature. In the objective space, diversity, representativeness, and the proximity of solutions obtained by the MOEA to the $PF$ set have long been used as parameters to measure an algorithm's efficiency. However, for decision-makers, the $PF$ set is often too large, making it a tedious and challenging task to select a particular solution. In some cases, decision-makers already have a predefined set of preferences, which restricts the search for solutions to a specific region of the $PF$ set. Despite respecting optimality, solutions outside this region are irrelevant or inappropriate for implementation. Therefore, it is beneficial to define a Region of Interest (ROI) in the objective space a priori and restrict the search for solutions to this region. This procedure optimizes computational effort and produces better solutions in terms of convergence \cite{ivan1}. These characteristics are even more relevant when the dimension $m$ of the objective space exceeds three.

Many preference-based solvers have been proposed in the literature, designed to converge to a subset of Pareto-optimal solutions located in a given Region of Interest (ROI) of the Pareto Front defined a priori by the decision maker \cite{gecco24,YADAV2025,Bechikh2015}. With this new approach, one can avoid the main disadvantages of a posteriori methods. Defining the region of interest may be easier and more efficient for the decision-maker than modeling preferences based on specific parameters of a parameterized single-objective optimization problem.

\begin{figure}[!ht]
	\centering
	\includegraphics[width=.45\linewidth]{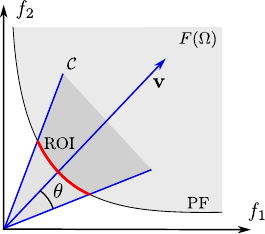}
	\caption{Preference cone $\mathcal{C}$ with its axis $\mathbf{v}$ and opening angle $\theta$, objective space $F(\Omega)$, Pareto Front ($PF$) and ROI in bi-dimensional space.}\label{fig:roi}
\end{figure}

In \cite{Bechikh2015}, the authors present several ways of defining an ROI, and MOEA applied to this context. \cite{ivan1} proposes the definition of an Interest Region using a preference cone $\mathcal{C}$ in the objective space. This cone is defined by a vector axis $\mathbf{v}$ and an opening angle $\theta$. On the $\mathbf{v}$-axis, the optimizer presents some relationship between its preferences in the form of a vector. For example, in a bi-objective minimization problem, the vector $\mathbf{v} = (1,1)$ defines a cone of the PF where the two objectives are equally important. In contrast, the vector $\mathbf{v} = (\sfrac{1}{2},1)$ defines a cone of the PF with preference for the first objective over the second objective\footnote{For maximization problems use $\mathbf{v} = (2,1)$}. The opening angle $\theta$ presents the aperture of the preference cone $\mathcal{C}$. The ROI is defined as $\mathcal{C} \cap PF$. Figure \ref{fig:roi} shows the ROI obtained by the intersection between cone $\mathcal{C}$ and PF in the two-dimensional objective space $F(\Omega)$. Small values of $\theta$ define a restricted region around the $\mathbf{v}$-axis. Increasing the value of $\theta$ expands the region of interest. Furthermore, concentrating solutions in an ROI in the objective space improves convergence indicators in problems with many objectives.

In \cite{gecco22}, the authors present a variation of the NSGAII algorithm to find solutions in a ROI. Their proposal defines the ROI using the minimum $f_i^L$ and maximum $f_i^U$ values of each objective $f_i$ of the application $F(\mathbf{x})$. Using these values, a hypercube $\mathcal{H}=[f_i^L \times f_i^U]$ is delimited in the objective space. Then, the ROI is defined as $\mathcal{H}\cap PF$.

Additionally, \cite{YADAV2025} proposed a framework for decision-making that selects preferred Pareto regions, rather than identifying the entire Pareto set, and addresses uncertainty in variables and problem parameters differently from the approach presented here.

\subsection{DWU algorithm}

The Dominance-Weighted Uniformity (DWU) approach, proposed in \cite{moreira:19}, is a multi-objective evolutionary algorithm designed to enhance solution diversity in the decision space. 
The uniformity measure computes the distance between the two closest points of a given set $R\subseteq \Omega$, and is formulated as
\begin{equation}\label{eq:U_R}
U(R)=\min_{\underset{r_1\neq r_2}{r_1,r_2\in R}} \| r_1 -r_2\| 
\end{equation}
where $\| \cdot \|$ is $L_2$ norm-based distance. Finding the Pareto optimal set with maximum dispersion is to solve the following problem
$$U(R^*)=\max_{\underset{|R|=k}{R\subseteq PO}} U(R)$$
given a set $PO$ of approximate efficient solutions generated by a multi-objective evolutionary algorithm, a predefined cardinality $k\leq |PO|$.

It incorporates a guiding criterion that strikes a balance between uniformity in the decision space and dominance information in the objective space. 
In the following, we describe the algorithm's structure.
\begin{enumerate}
\item Generate, randomly, the initial population $\mathcal P_0$ of size $N$.
\item Non-dominated sorting levels $FL_1, FL_2, \ldots FL_D$ are computed for all the $N$ individuals of $P_0$.
\item While $t \leq max\_generation$, the population $\mathcal P_{t+1}$ is obtained as follows: 
\begin{itemize}
\item From population $\mathcal P_t$, perform $k$ 
binary tournaments, crossover and mutation operators, 
generating a list $\mathcal L_t$ of $N$ new individuals.
\item A combined population $\mathcal A_t = \mathcal P_t \cup \mathcal L_t$ of size $2N$ is obtained.
\item From set $\mathcal A_t$, 
the population $\mathcal P_{t+1}$,
with cardinality $N$, it is selected by the DWU heuristic.
\item $t \longleftarrow t + 1$
\end{itemize}
\item Return final non-dominated solutions.
\end{enumerate}

The selection process employs a greedy heuristic, the DWU heuristic (Algorithm \eqref{alg:mod}), which prioritizes solutions with high uniformity while considering their dominance level, as given in the MOP (Equation \eqref{eq:mop1}).
\begin{algorithm}[!th]
\caption{DWU heuristic}
\label{alg:mod}
\begin{algorithmic}[1]
\STATE \textbf{input}: \textit{Feasible set $\mathcal X\subset \Omega$ of the MOP (\ref{eq:mop1}) and scalar $k$.}
\STATE $N_\mathcal{X}\subseteq \mathcal{X}$ such as $F(N_\mathcal{X})$ is the set of non-dominated solutions.
\STATE Let $x_i, x_j$ in $N_\mathcal X$ with maximum $w_d( x_i, x_j)$.
\STATE $R \longleftarrow \{x_i,x_j\}$
\WHILE{$|R|<k$}
\STATE Find $x^\prime\in \mathcal{X} \setminus R$ such that $\displaystyle \min_{r\in R}{w_d( x^\prime, r)}$ is maximum.
\STATE $R \longleftarrow R \cup \{x^\prime\}$
\ENDWHILE
\STATE \textbf{output}: $R$
\end{algorithmic}
\end{algorithm}

The authors define the $w_d()$, dominance-weighted uniformity function, to select solutions uniformly spread in the decision space. Given $x, x^{''}\in \mathcal{X}$:
\begin{equation}
w_d(x,x^{''})=\frac{\| x-x^{''}\|}{|d(x)-d(x^{''})|+1}
\end{equation}
wherein 
\begin{equation}\label{spea}
    d(x)=\sum_{
    \substack{x^{''} \in \mathcal \mathcal{X}\\
    \mathbf{F}(\mathbf{x^{''}})\prec \mathbf{F}(\mathbf{x})}}{s(x^{''})}
\end{equation}
wherein $s(x)$ in the population ${\mathcal X}$ is defined by the number of solutions in $\mathcal X$ it dominates and, $\| \cdot \|$ is $L_2$ norm-based distance (Euclidean norm) take up in the decision space. 

\section{Proposed Approach}

\begin{figure}[!ht] 
	\centering
	\includegraphics[width=.4\textwidth]{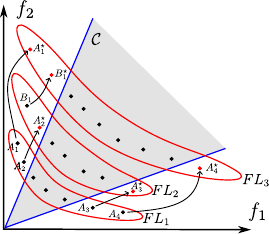}
	\caption{Penalization of solutions located outside the preference cone $\mathcal{C}$} \label{fig:NSGAII-penal}
\end{figure}

We propose the $\mathcal{C}$-DWU to incorporate the preference cone information supplied by the decision-maker into the DWU algorithm ~\cite{moreira:19}.
The $\mathcal{C}$-DWU requires only the penalized threshold values of the solutions located outside the preference cone $\mathcal{C}$ with its axis $\mathbf{v}$ and opening angle $\theta$ modified two steps of the DWU algorithm proposed in \cite{moreira:19}, as follows:

\noindent
\textbf{Non-dominated sort levels:} 
Figure \ref{fig:NSGAII-penal} illustrates the penalization process idea for solutions located outside the preference cone $\mathcal{C}$, in the non-dominated sorting phase. The solutions $A_1$, $A_2$, $A_3$ and $A_4$, located in the first non-dominated front level ($FL_1$), should be penalized. Solutions $A_2$ and $A_3$ are reclassified in the second non-dominated front level ($FL_2$) because they are at a smaller angular distance than solutions $A_1$ and $A_4$, which are reclassified to non-dominated front level 3 ($FL_3$). These reclassified solutions are denoted by $A_1^{\star}$, $A_2^{\star}$, $A_3^{\star}$ and $A_4^{\star}$ respectively. The same situation occurs with solution $B_1$, located on the second non-dominated front level ($FL_2$), reclassified to the third non-dominated front level ($FL_3$) and denoted by $B_1^{\star}$. 

The non-dominated sorting phase of $\mathcal{C}$-DWU penalizes the solutions located outside the preference cone, proportional to the angular distance between $\phi$, the solution, and the cone axis $\mathbf{v}$, as follows:
\begin{equation}\label{eq:penal_front}
	\mathcal{C}\mbox{-}FrontLevel = \begin{cases}
    \begin{tabular}{@{}cc@{}}
        $FrontLevel + P_{\alpha, \theta}(\phi) $,
        & $\mbox{if } \phi > \theta$ \\
        &\\
         $FrontLevel$, & $\mbox{if } \phi \leq \theta $
    \end{tabular}\\
	\end{cases}
\end{equation}
wherein $P_{\alpha, \theta}(\phi) =
    \left \lfloor{ 
    e^{\left(\alpha(\phi - \theta)\right)}}
    \right \rfloor$
is the penalization function ($\lfloor \star  \rfloor$ is the floor function (rounds to the nearest integer less than or equal to)), with  $\theta$ (opening angle of the preference cone $\mathcal{C}$) and $\alpha$ (penalty intensity) parameters. $FrontLevel$ refers to the level of non-domination of the solution within the population. 

\noindent
\textbf{Dominance-weighted uniformity function:} Similarly, the idea is to penalize solutions outside the preference cone $\mathcal{C}$ in the selection phase, which considers the uniformity of solutions in the decision space.

We propose a penalty value to the dominance-weighted uniformity function proposed in~\cite{moreira:19}, named $\mathcal{C}\mbox{-}w_d$, as follows:
\begin{equation}\label{eq:penal_dwu}
\mathcal{C}\mbox{-}w_d(x,x^\prime)=\begin{cases}
    \begin{tabular}{@{}c@{}c@{}}
        $\displaystyle\frac{\| x-x^\prime\|}{|d(x)-d(x^\prime)|+1}- P_{\beta, \theta}(\phi)$,& $\mbox{if } x \mbox{ or }x^\prime \notin \mathcal{C}$ \\
        &\\
         $\displaystyle\frac{\| x-x^\prime\|}{|d(x)-d(x^\prime)|+1}$, & $\mbox{if }x \mbox{ and }x^\prime \in \mathcal{C}$
    \end{tabular}\\
\end{cases}
\end{equation}
wherein $\mathcal{C}$ is the preference cone defined by axis $\mathbf{v}$ and angle $\theta$, $P_{\beta, \theta}(\phi)=e^{(\beta(\phi - \theta))}$, where $\phi$ is the angular distance between the image of the solutions $\mathbf{x}, \mathbf{x}^\prime$ and the cone axis $\mathbf{v}$. The parameter $\beta $ determines the penalty intensity.  

\begin{figure}[!ht]
    \begin{subfigure}{0.48\textwidth}
    \centering
        \includegraphics[width=.9\linewidth]{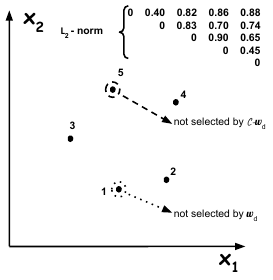}
        \caption{Hypothetical example with five solutions in the decision space and respective distances.}
    \end{subfigure}\hfill
    \begin{subfigure}{0.48\textwidth}
    \centering
	\includegraphics[width=.9\linewidth]{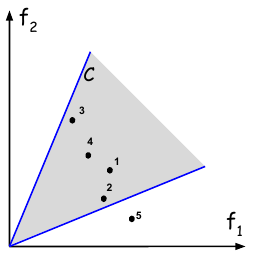}
    \caption{Images of the solutions in the objective space, for a bi-objective minimization problem.}
    \end{subfigure}
    \caption{Symbolic application of the proposed
$\mathcal{C}\mbox{-}w_d$. An example to find a representative subset $R$, with $|R| = 4$, of a set $\mathcal{P}$, where $|\mathcal{P}| = 5$. Figure (a) shows the solution
in the decision space not selected by $\mathcal{C}\mbox{-}w_d$ and $w_d$.}\label{fig:cdw-penal}
\end{figure}

The idea is to integrate the information of the solution in objective space into the diversity indicator in decision space. To this end, we consider the example of finding a representative subset $R$, with $|R| = 4$, of a set $\mathcal{P}$, where $|\mathcal{P}| = 5$, shown in Figure \ref{fig:cdw-penal}. The penalization proposed in equation \ref{eq:penal_dwu} takes up the information about the image of the solution 5, which is discarded in the selection process.

The goal is to find the solutions that ``maximize'' the concentration in the objective space and also ``maximize'' the dispersion in the decision space, ensuring $\mathbf{F}(\mathbf{x})$ in the Pareto-optimality context.

\section{Experiments}

In this study, we define a bi-objective constrained optimization problem (BOP), as follows:
\begin{equation}
	\begin{array}{l}
		\min ~~ \mathbf{F}(\mathbf{x})=\left(f_1(\mathbf{x}),f_2(\mathbf{x})\right)\\
		\\
		\mbox{subject to:} ~~ 
		\begin{cases}
			\displaystyle\arccos{\left
				(\frac{\left<{\mathbf{F}}(\mathbf{x}),{v}\right>}{\|{\mathbf{F}}(\mathbf{x})\| \|{v}\|}\right)}\leq \theta \\
			\mathbf{x}=(x_1,x_2, \ldots, x_n) \in \Omega \subset \mathbb{R}^n
		\end{cases}
	\end{array} 
	\label{eq:multi1}
\end{equation}

\noindent
\textbf{Maximum Concentration of FP}: The constraints defined in the \eqref{eq:multi1} define the direction of search over the preference cone $\mathcal{C}$ in the objective space, named the concentration in the objective space.

\noindent
\textbf{Maximum Dispersion of PO:} We used the direction of search in the decision space as the dispersion measure presented in \cite{moreira:19}, given by equation (\ref{eq:U_R}).

Considering two objectives, we applied the proposed methodology to multimodal problems (WFG4 and WFG9) and an unimodal problem (DTLZ2). In each problem, we set the decision space dimension to 5, 7, and 9. The population size for each algorithm was set to 100 individuals, with $10^5$ objective function evaluations. The preference cone $\mathcal{C}$ in the objective space was defined with an axis $\mathbf{v}=(1,1)$ and an opening angle of $\theta = 0.3$ radians. The parameters $\alpha$ and $\beta$ (see Equations \eqref{eq:penal_front} and \eqref{eq:penal_dwu}) of the proposed algorithm $\mathcal{C}$-DWU were empirically defined as $\alpha = 0.3$ and $\beta= 1$ determines the penalty intensity. All experiments were conducted using the PLATEMO platform \cite{Tian2017}. The code and results are available in \url{https://github.com/moreirag/MDMC-MOP}.

We analyzed the adequacy of the solutions obtained in the Region of Interest (ROI) based on whether they were included or excluded within the preference cone. The convergence of each algorithm was measured using the Inverted Generational Distance (IGD) metric, considering the ROI Pareto Front, i.e., part of the Pareto Front contained in the ROI. The uniformity metric proposed in \cite{moreira:19} was employed to assess the dispersion of solutions in the Pareto Optimal set ($PO$).

We compared the proposed algorithm against a variation of NSGA-II~\cite{Deb2002}. Similarly, the penalization (\ref{eq:penal_front}) is added to the dominance front $FrontLevel$ of the solutions after the non-dominated sorting mechanism of NSGA-II. The standard crowding distance function and the algorithm proceed as usual. Let us denote this variation by $\mathcal{C}$-NSGAII.

\section{Results}

We conducted 10 runs for each algorithm, instance, and problem combination. Across all simulations, the constraint ensuring that the solutions belong to the ROI was violated in only one instance of $\mathcal{C}$-DWU/WFG9, with $D = 9$. The mean IGD values are presented in Table \ref{tbl:igd}, and Table \ref{tbl:u} presents the mean Uniformity Measure values.

Table \ref{tbl:igd} shows that both methodologies present good convergence in the $PF$ set, with superior results favoring the $\mathcal{C}$-NSGAII algorithm.
\begin{table}[!ht]
	\centering
	\caption{Mean IGD value in 10 runs of the $\mathcal{C}$-NSGAII and $\mathcal{C}$-DWU algorithms.}
    \renewcommand{\tabcolsep}{8 pt}
		\begin{tabular}{@{}lcccc@{}}
		\toprule
		\multirow{2}{*}{Algorithm} & \multirow{2}{*}{Problem} & \multicolumn{3}{c}{Decision Space Dimension (D)} \\ \cmidrule(l){3-5} 
		&  & 5 & 7 & 9 \\ \midrule
		\multirow{3}{*}{$\mathcal{C}$-DWU} & DTLZ2 & 2.9807e-03 & 3.2606e-03 & 3.3158e-03 \\
		& WFG4 & 9.4183e-03 & 9.4628e-03 & 1.0037e-02 \\
		& WFG9 & 1.0169e-02 & 1.6987e-02 & 1.8381e-02 \\ \midrule
		\multirow{3}{*}{$\mathcal{C}$-NSGAII} & DTLZ2 & 1.8003e-03 & 1.9312e-03 & 1.9006e-03 \\
		& WFG4 & 5.7448e-03 & 5.6247e-03 & 5.8095e-03 \\
		& WFG9 & 6.8577e-03 & 9.4836e-03 & 8.7864e-03 \\ \bottomrule
	\end{tabular}
    \label{tbl:igd}
\end{table}

However, analyzing Table \ref{tbl:u}, it can be observed that the $\mathcal{C}$-DWU algorithm achieved better results regarding the average uniformity of solutions in the $PO$ set.
\begin{table}[!ht]
\centering
	\caption{Mean uniformity measure value in 10 runs of the $\mathcal{C}$-NSGAII and $\mathcal{C}$-DWU algorithms.}
	\label{tbl:u}
    \renewcommand{\tabcolsep}{8 pt}
	\begin{tabular}{@{}lcccc@{}}
		\toprule
		\multirow{2}{*}{Algorithm} & \multirow{2}{*}{Problem} & \multicolumn{3}{c}{Decision Space Dimension (D)} \\ \cmidrule(l){3-5} 
		&  & 5 & 7 & 9 \\ \midrule
		\multirow{3}{*}{$\mathcal{C}$-DWU} & DTLZ2 & 0.0251 & 0.0299 & 0.0338 \\
		& WFG4 & 0.1081 & 0.2466 & 0.4045 \\
		& WFG9 & 0.0949 & 0.3145 & 0.4416 \\ \midrule 
		\multirow{3}{*}{$\mathcal{C}$-NSGAII} & DTLZ2 & 4.4293e-04 & 3.9101e-04 & 7.0981e-04 \\
		& WFG4 & 3.7817e-05 & 3.6023e-05 & 7.5732e-05 \\
		& WFG9 & 4.3458e-04 & 2.7848e-04 & 7.1390e-04 \\ \bottomrule
	\end{tabular}
\end{table}

\begin{figure}[!ht]
    \begin{subfigure}[b]{0.48\textwidth}
    \centering
        \includegraphics[width=\linewidth]{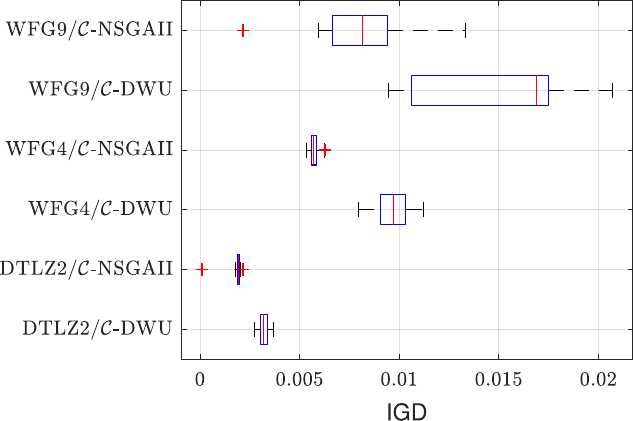}
        \caption{IGD value of the $\mathcal{C}$-NSGAII and $\mathcal{C}$-DWU algorithms.}
	\label{fig:IGD}
    \end{subfigure}
	\hfill
    \begin{subfigure}[b]{0.48\textwidth}
	\centering
	\includegraphics[width=\linewidth]{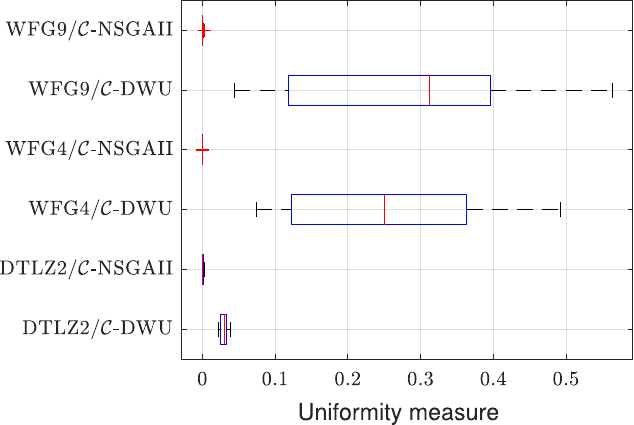}
    \caption{Uniformity measure value of the $\mathcal{C}$-NSGAII and $\mathcal{C}$-DWU algorithms.}
	\label{fig:U}
    \end{subfigure}
    \caption{Box plot of the IDG and uniformity measure value in 10 runs of the $\mathcal{C}$-NSGAII and $\mathcal{C}$-DWU algorithms.}
\end{figure}

Figures \ref{fig:IGD} and \ref{fig:U} show the boxplots of the metrics' values. In Figure \ref{fig:IGD}, we can compare the values of the IGD metric obtained in the experiments. The $\mathcal{C}$-NSGAII algorithm obtained better results, confirming what was presented in Table \ref{tbl:igd}.

Figure \ref{fig:U} shows the measurements of the uniformity metric applied to the solutions of the set $PO$. A significant difference is observed in the metric values for the multimodal functions (WFG4 and WFG9). The solutions obtained by the $\mathcal{C}$-NSGAII algorithm present lower values for the uniformity metric and little variability in these values. On the other hand, the solutions obtained by the $\mathcal{C}$-DWU algorithm, in addition to presenting higher values of the dispersion measure, also show greater variability, indicating the greater dispersion of the solutions in the set $PO$. Even in the unimodal problem DTLZ2, the dispersion measures obtained by the $\mathcal{C}$-DWU algorithm exceed the values found by the $\mathcal{C}$-NSGAII algorithm. The lower variability of the dispersion measures for the solutions to this problem is explained by its unimodal nature, which means that a single region of the set $PO$ concentrates on the optimal solutions. Despite this, the proposed methodology could still find dispersed solutions, preserving the convergence of solutions in the objective space.

\begin{figure}[ht!]
    \begin{subfigure}[b]{0.48\textwidth}
    \centering
        \includegraphics[width=.9\textwidth]{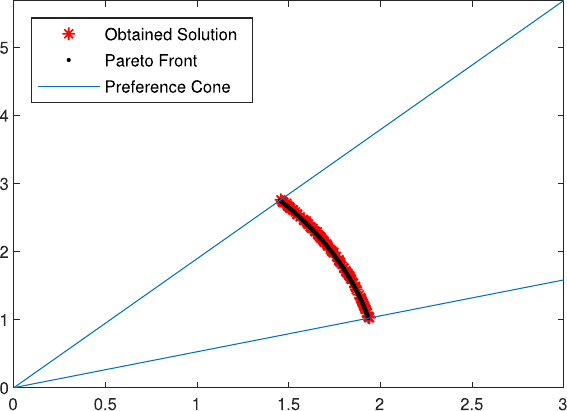}	
    \caption{Objective space for the WFG4 Problem using the $\mathcal{C}$-DWU algorithm.}
	\label{fig:dwu-obj}
    \end{subfigure}
	\hfill
    \begin{subfigure}[b]{0.48\textwidth}
	\centering
	\includegraphics[width=.9\textwidth]{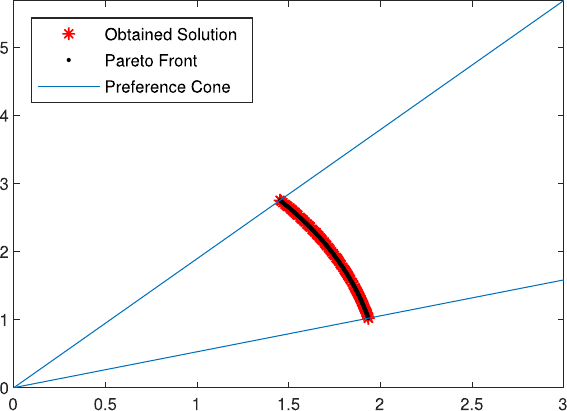}
    \caption{Objective space for the WFG4 Problem using the $\mathcal{C}$-NSGAII algorithm.}
	\label{fig:nsgaii-obj}
    \end{subfigure}
    \caption{Solutions obtained in the objective space for the two-objective problem WFG4.}
\end{figure}

Figures \ref{fig:dwu-obj} and \ref{fig:nsgaii-obj} illustrate the solutions obtained in the $PF$ set for the WFG4 function with 9 (nine) decision variables using the $\mathcal{C}$-DWU and $\mathcal{C}$-NSGAII algorithms, respectively. It can be observed that the difference in the IGD metric between these two algorithms does not significantly impact the proximity between the solutions obtained and the Pareto Front.

\begin{figure}[!ht]	
    \begin{center}
    \begin{subfigure}{0.48\textwidth}
    \centering
    \includegraphics[width=.8\textwidth]{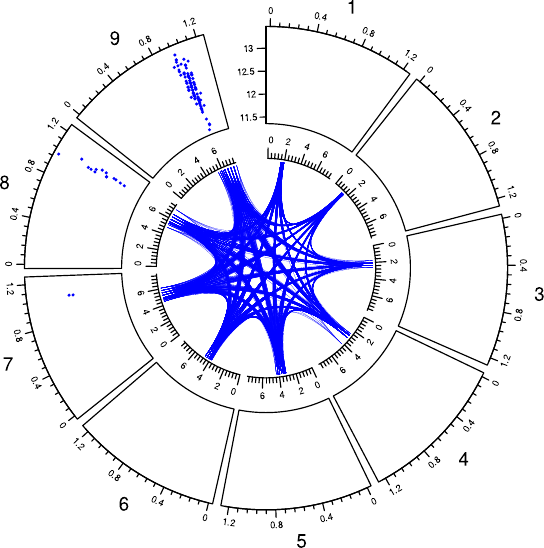}		
    \caption{Spatial dispersion visualization of solutions in the decision space for the WFG4 Problem using the $\mathcal{C}$-DWU algorithm.}
	\label{fig:dwu-dec-cap}
    \end{subfigure}
	\hfill
    \begin{subfigure}{0.48\textwidth}
    \centering
	\includegraphics[width=.8\textwidth]{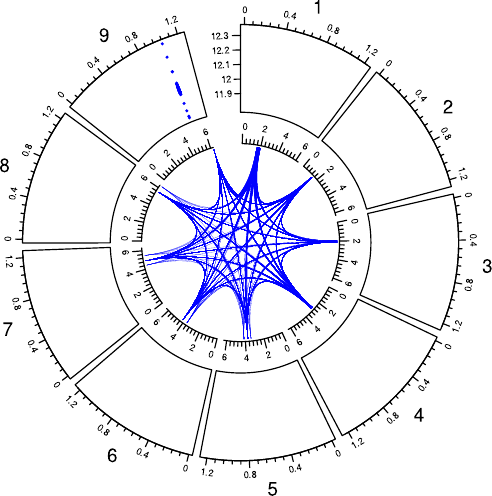}	
    \caption{Spatial dispersion visualization of solutions in the decision space for the WFG4 Problem using the $\mathcal{C}$-NSGAII algorithm.}
	\label{fig:nsgaii-dec-cap}
    \end{subfigure}
    \end{center}
    \vspace{1em}
    \begin{subfigure}{0.48\textwidth}
    \centering
    \includegraphics[width=\textwidth]{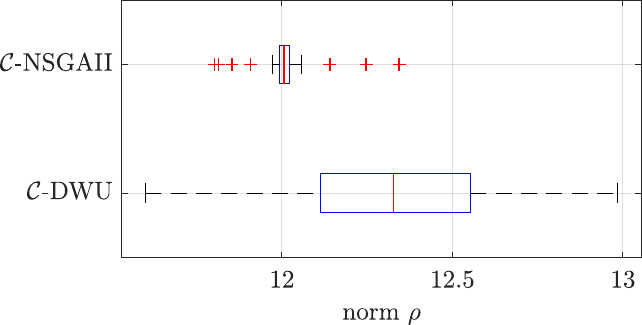}
    \caption{Box plot of the norm $\rho$ of solutions in the decision space for the WFG4 problem using the $\mathcal{C}$-DWU algorithm.}
	\label{fig:boxplot_norma}
    \end{subfigure}
	\hfill
    \begin{subfigure}{.48\textwidth}
	\centering
	\includegraphics[width=\textwidth]{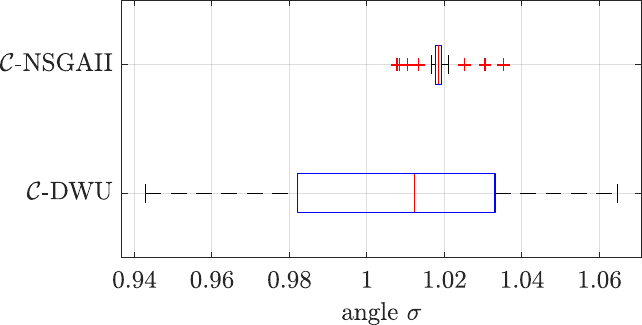}
    \caption{Box plot of the smallest angle $\sigma$ of solutions in the decision space for the WFG4 problem using the $\mathcal{C}$-NSGAII algorithm.}
	\label{fig:boxplot_angulo}
    \end{subfigure}
    \caption{Spatial dispersion visualization of solutions obtained in the decision space for the two-objective problem WFG4.}
    \label{fig:boxplot_norma_angulo}
\end{figure}

Figures \ref{fig:dwu-dec-cap} and \ref{fig:nsgaii-dec-cap} show the spatial dispersion visualization of solutions in the decision space. These figures utilize the CAP-VIZ tool for high-dimensional data visualization \cite{Koochaksaraei2017,8477889}. In this methodology, the points of the decision space are classified according to the smallest angular distance $\sigma$ relative to each of the axes $e_1=(1,0, \ldots, 0), e_2=(0,1,0, \ldots,0), \cdots, e_9=(0, \ldots, 0,1)$ in the decision space $\mathbb{R}^9$. After classification, the Euclidean norm $\rho$ of each point is computed. The pair $(\sigma, \rho)$ is recorded in the outer section of the plot corresponding to the axis where the smallest angle $\sigma$ was observed. Each of these 9 (nine) regions is referred to as a sector of the plot. In the outer portion of each sector, the radial scale measures the norm $\rho$, while the scale surrounding the sector measures the angle $\sigma$. The largest value on this scale corresponds to the angle between an axis $e_i$ and the diagonal $d=(1,1, \ldots, 1)$ of the hypercube defined by the axes $e_1, \ldots, e_9$. The central region illustrates the chord diagram of all points.

The comparison between Figures \ref{fig:dwu-dec-cap} and \ref{fig:nsgaii-dec-cap} reveals that the solutions obtained in the decision space by the $\mathcal{C}$-NSGAII algorithm (Figure \ref{fig:nsgaii-dec-cap}) are concentrated solely in sector 9, meaning they are closer to the $e_9=(0,0, \ldots,1)$ axis of this space. The solutions form a cluster of aligned solutions, with some isolated points. Their spatial arrangement presents little diversity, and consequently, the decision-maker is presented with few options.

On the other hand, the solutions obtained using the $\mathcal{C}$-DWU algorithm, illustrated in Figure \ref{fig:dwu-dec-cap}, are distributed across sectors 7, 8, and 9, with a more substantial presence in the latter. The solutions obtained within sector 9 are more dispersed than those in the same sector in Figure \ref{fig:nsgaii-dec-cap}. No clusters or isolated points are observed, except for solutions in sector 7. In addition to the greater variability of solutions present in sector 9, solutions are also observed in sector 8, which corresponds to the region of space closest to the axis $e_8=(0,0, \ldots, 1,0)$, two isolated solutions in sector 7 (region near the axis $e_7=(0, \ldots, 1,0,0)$). This situation presents more options to the decision maker, even allowing a fine-tuning of the chosen solutions. The discrepancy between the variability of the $\rho$ norm and the smallest angle $\sigma$ between the two proposals is presented in Figure \ref{fig:boxplot_norma} and \ref{fig:boxplot_angulo} respectively.

The variability of the obtained solutions is also evident in the chord diagram, positioned at the center of the graphic. In Figure \ref{fig:nsgaii-dec-cap}, the ends of the arcs are clustered within a few coordinates. In sectors 2, 3, 4, and 9 (corresponding to the same coordinates as the decision space solutions), the vertices of the arcs originate from a single measurement neighborhood. In sectors 5, 6, 7, and 8, the ends of the arcs stem from the neighborhood of two measurements (in sectors 6 and 8) or three measurements (in sectors 5 and 7).

In contrast, the chord diagram in Figure \ref{fig:dwu-dec-cap} presents arc vertices distributed over wider intervals. In particular, in sectors 7, 8, and 9, a greater dispersion of the values of these coordinates is observed when compared to the same coordinates of the solutions obtained by the $\mathcal{C}$-NSGA-II algorithm. This indicates that in each parameter of the decision space, the range of choice of solutions will be larger, demonstrating greater flexibility in selecting some of the optimal solutions obtained.

\section{Conclusion}

This paper presents the $\mathcal{C}$-DWU solver for constrained bi-objective optimization problems, which simultaneously concentrates in a Region of Interest in the objective space while exhibiting maximum dispersion in the decision space. We also present a solution methodology and analysis of the obtained results. The proposed problem illustrates a common real-world scenario where a decision-maker requires solutions with maximum diversity in their parameters and has sufficient understanding of the problem to define a Region of Interest in the objective space.

The proposed methodology successfully obtains high-quality solutions in the objective space within the predefined Region of Interest with good convergence. It also ensured that the optimal solutions exhibited significantly higher diversity than those obtained by the modified NSGAII algorithm. Defining the Region of Interest(ROI)  using the preference cone $\mathcal{C}$ is an available alternative to other methods due to its ease of adjusting the direction and extent of the ROI, as well as its straightforward incorporation into existing algorithms.

The trade-off between the convergence of solutions in the set $PF$ and the dispersion in the set $PO$ in the obtained solutions presents a better situation for the dispersion of solutions. The reduction in convergence rates is much lower than the gain in solution diversity. Diverse solutions help mitigate the bias caused by clustering in a specific region of the $PO$ set. In addition, it allows the decision maker to choose new optimal solutions when faced with situations not foreseen before the computational modeling of the problem, which may make specific regions of the set $PO$ rich in solutions unfeasible.

\section*{Funding}
The authors would like to thank the Fundação de Amparo a Pesquisa do Estado de Minas Gerais (FAPEMIG, grant APQ-01647-22), Conselho Nacional de Desenvolvimento Científico e Tecnológico (CNPq, grants 307151/2022-0, 152613/2024-2) and Instituto Federal de Educação e Tecnologia de Minas Gerais (IFMG, grant 030/2024) for supporting the development of this study.





\bibliographystyle{main}  
\bibliography{main}

\end{document}